\documentclass[11pt]{article}

\usepackage{graphicx} % Required for inserting images
\usepackage{amsfonts}
\usepackage{amsmath}
\usepackage{amsthm}
\usepackage{color}
\usepackage{hyperref}
\newtheorem{definition}{Definition}[section]
\newtheorem{theorem}{Theorem}[section]
\newtheorem{lemma}{Lemma}[section]
\newtheorem{assumption}{Assumption}[section]
\newtheorem{remark}{Remark}[section]

\title{Extended HJB Equation for Mean-Variance Stopping Problem: Vanishing Regularization Method}
\author{Yuchao Dong\thanks{{ School of Mathematical Sciences, Key Laboratory of Intelligent Computing and Applications (Ministry of Education), Tongji University, Shanghai 200092, China. The work of the second author was funded by National Natural Science Foundation of China
			(No.12471425) }} 
and Harry Zheng\thanks{Department of Mathematics, Imperial College, London SW7 2BZ, UK (Email: h.zheng@imperial.ac.uk). The work was supported by   the Engineering and Physical Sciences Research Council of
UK  (Grant No. EP/V008331/1).}}
\begin{document}

\maketitle
\abstract{
This paper studies the time-inconsistent MV optimal stopping problem via a game-theoretic approach to find equilibrium strategies. To overcome the mathematical intractability of direct equilibrium analysis, we propose a vanishing regularization method: first, we introduce an entropy-based regularization term   to the MV objective, modeling mixed-strategy stopping times using the intensity of a Cox process. For this regularized problem, we derive a coupled extended Hamilton-Jacobi-Bellman (HJB) equation system, prove a verification theorem linking its solutions to equilibrium intensities, and establish the existence of classical solutions for small time horizons via a contraction mapping argument. By letting the regularization term tend to zero, we formally recover a system of parabolic variational inequalities that characterizes equilibrium stopping times for the original MV problem.   This system includes an additional key quadratic term--a distinction from classical optimal stopping, where stopping conditions depend only on comparing the value function to the instantaneous reward.

}



\section{Introduction}
Given a diffusion process $X$, the classical optimal stopping problem is to determine a stopping time $\tau$ that maximizes 
$$
\mathbb E\left[f(X_{\tau})\right]. 
$$ 
Optimal stopping has many applications, for example, financial decision-making (e.g., timing for asset sales) and statistical inference (e.g., stopping rules for hypothesis testing \cite{tartakovsky2014sequential}). However, in financial contexts, there is often an additional imperative to mitigate decision-related risk. In line with the mean-variance analysis of \cite{markowitz2008portfolio}, we identify the return with the expectation  and the risk with the variance and aim to select a stopping time that maximizes
\begin{equation} \label{MVstopping}
\mathbb E\left[f(X_{\tau})\right]-\frac\gamma 2\text{Var}\left[f(X_{\tau})\right],
\end{equation}
where $\gamma \ge 0$ denotes the risk aversion coefficient. Problem (\ref{MVstopping})  is referred as a mean-variance (MV) stopping problem in the literature. 

Similar to the classical dynamic MV problem, the optimal stopping rule  typically depends on the initial state $x$, which means that it does not generally satisfy Bellman’s principle of optimality. In the literature this is known as {\it time-inconsistency}. Time-inconsistent problems are typically studied using two  approaches. One is to formulate the problem for a fixed initial state and allow the corresponding optimal stopping rule to depend on that initial state, called the pre-commitment approach. The other is to look for a stopping rule that remains optimal at every period when re-evaluated from that period's perspective, called the game-theoretic approach.

Strotz \cite{strotz1955myopia} is the first to explore the game-theoretic approach to time-inconsistent problems in dynamic utility maximization with non-exponential discounting. Bjork et al. \cite{bjork2021time} give a comprehensive treatment of time-inconsistent Markovian models   and characterize the  equilibrium by a solution to a generalized HJB equation, called the extended HJB system.  Time inconsistent control problems have attracted considerable research interest  in recent years with many applications. For example,  Bjork et al. \cite{bjork2014mean} solve a MV problem with state dependent risk aversion. He and Liang \cite{he2013optimal} study a defined contribution insurance problem in a MV framework. Dai et al. \cite{dai2023learning} solve a MV problem with reinforcement learning method. 
All aforementioned papers have fixed finite horizon.

The literature on the game-theoretic approach to time-inconsistent stopping problems is in the early developing stage. Christensen and Lindensj\"o \cite{christensen2018finding}  study an equilibrium stopping problem with initial state dependent   reward.  Bayraktar et al. \cite{bayraktar2023equilibria} consider three equilibrium concepts proposed in  the literature for time-inconsistent stopping problems with non-exponential discount. 
There is little research   for MV stopping problems.  The only ones the authors are aware of are Peskir and Shiryaev  \cite{peskir2006optimal} on the so-called dynamic optimal stopping time, which is similar to the game theoretic approach and  Christensen and Lindensj\"o \cite{christensen2020time} on a subgame perfect Nash equilibrium for stopping problems.

In this paper, we study the equilibrium strategy and relate it to  the extended HJB equation, which means that we need  to formulate the problem as a game and look for equilibrium.   It is a fundamental result in game theory that equilibrium generally exists for mixed strategies rather than pure strategies in a broad class of games\footnote{For example, the rock - paper - scissors game is a classic example in game theory, and it has no pure strategy Nash equilibrium but has a mixed  strategy Nash equilibrium in which  the player choose each action with equal probability $1/3$. }. Hence,  we focus on mixed strategy stopping times by allowing the agents to  choose the intensity function of a Cox process as a randomization device for the stopping decision, whereas \cite{christensen2020time, peskir2006optimal} characterize the equilibrium and  provide other necessary and sufficient equilibrium conditions, but do not derive the extended HJB equation. While their results coincide with ours for geometric Brownian motion case, our derivation is motivated by the vanishing regularity approach, a key distinction from prior research.

We now describe the key methodology for solving the MV stopping problem.  We first add a regularization term, weighted by a constant $\lambda$ into the target functional (see (\ref{J^lambda})), to form a regularized problem with control variate being the intensity (as opposed to the stopping time), which makes   the definition of the equilibrium  straightforward. We then derive the associated extended HJB equation (see (\ref{extended-HJB})), prove a verification theorem (see Theorem \ref{thm_equil})  and establish the existence of its solution under certain technical assumptions (see Theorem \ref{thm3.2}). Finally, we  let $\lambda$ tend to zero (i.e., vanishing regularization) to formally obtain a system of parabolic variational inequalities (see (\ref{extended_VI})) that characterizes the equilibrium stopping time for the original MV stopping problem. To the best knowledge of the authors, this is the first time  such a system of equations is reported  in the literature. Furthermore, we find that the stopping condition is not determined solely by comparing the value function with the instantaneous reward; instead, an additional quadratic term also plays a role in the formulation, which is in sharp contrast to the standard  optimal stopping problem.

Finally, we aim to emphasize the motivation underlying our method. In \cite{dong2024randomized}, the authors investigated an entropy-regularized optimal stopping problem and demonstrated that the corresponding optimal value function is associated with a penalized form for the variational inequality. Notably, this penalized equation also converges to the original variational inequality as the regularization parameter tends to zero.	For MV optimal stopping problems, we adopt an analogous approach: we first introduce and analyze a regularized version of the problem, then subsequently let the regularization parameter vanish to recover results for the original MV stopping problem. Such ideas are pervasive in mathematical research. When addressing a computationally or theoretically challenging problem, researchers often first consider a perturbed or regularized counterpart—whose solution is more tractable to derive—and then take an appropriate limit to revert to the original problem. It is precisely this core logic that leads us to name our approach {\it Vanishing Regularization Method}.

The rest of the paper is organized as follows: In Section 2 we formulate the MV stopping problem and its relaxed problem together with the definition of equilibrium. In Section 3 we derive the extended HJB equation for the relaxed problem and prove a verification theorem (Theorem \ref{thm_equil}) and an existence result (Theorem \ref{thm3.2}). In Section 4 we let $\lambda$ tend to zero to formally get the variation system for the original MV stopping problem and show that it characterizes the equilibrium stopping time (Theorems \ref{thm_stopping} and \ref{thm_stopping1}). In Section 5 we give some further discussions on our results, including infinite horizon case, discrete time approximation, and general time-inconsistent problems. Section 6 concludes the paper. Appendix contains the proofs of a local time approximation relation (\ref{esti_square_of_local_time}) and a technical lemma  (Lemma \ref{append_lemma}) that is needed in the proof of Theorem \ref{thm3.2}.

\section{Mean Variance Stopping  and its Relaxed Problem}
In this section, we  introduce the basic framework of the problem. Especially,  we give the formulation for the relaxed MV stopping problem and the definition of related equilibrium strategy.  
Let $(\Omega,\mathcal F,P)$ be a probability space, on which a standard $d$-dimensional Brownian motion $W$ is defined\footnote{For simplicity, we consider the case that the dimension of the Brownian motion is same to that of the state process. It can be extended to other cases without any major modification as long as Assumption \ref{ass} holds.}. The $\mathbb R^d$-valued state process $X$ satisfies
$$
dX_t=b(t,X_t)dt+\sigma(t,X_t)dW_t.
$$
Denote by $\mathbb F:=\{\mathcal F_t\}_t$ the natural filtration generated by $W$, augmented by all $P$-null sets. The set $\mathcal T_{t,T}$ is defined as the totality of all $\mathbb F$-stopping time taking values in $[t,T]$. For any time $t \in[0,T]$, the  MV stopping problem is to choose $\tau \in \mathcal T_{t,T}$ such that the following  functional 
$$
J(t,x;\tau):=\mathbb E_{t,x}\left[f(X_\tau)\right]-\frac{\gamma}{2}\text{Var}_{t,x}\left[f(X_\tau)\right]
$$
is maximized, where $\mathbb E_{t,x}[\cdot]$ and $\text{Var}_{t,x}[\cdot]$ denote the conditional expectation and conditional variance conditioning on $X_t=x$, respectively.

In this paper, we shall adopt the following assumptions on the coefficients.
\begin{assumption}\label{ass}
The coefficients $b,\sigma$ and $f$ are   Lipschitz continuous  with linear growth  in  $x$, uniformly in $t$, i.e., there exists a constant $C$ such that, for any $t\in [0,T]$ and $x,y \in \mathbb R^d$, 
$$
|b(t,x)|,|\sigma(t,x)|,|f(x)| \le C(1+|x|),
$$
and
$$
|b(t,x)-b(t,y)|,|\sigma(t,x)-\sigma(t,y)|,  |f(x)-f(y)| \le C|x-y|.
$$
Moreover, $\sigma\sigma^T$ is uniformly non-degenerate, i.e., there exists a constant $c$ such that $\sigma\sigma^T(t,x) \ge cI$ for any $t\in [0,T]$ and $x \in \mathbb R^d$.  

\end{assumption}

For any function $f$ defined on $[0,T]\times \mathbb R^d$, we use the notations $\partial_t f$  and $\partial_x f$  to represent the derivatives with respect to $t$ and $x$ and $\partial_{xx} f$ its Hessian.  In this paper, we use $C$ to represent a constant  that could depend on the coefficients but  may be different from line to line.

It is well-known that, for  MV problem, the dynamic programming principle fails. Thus, people focus on two kinds of strategies. One is called pre-committed  strategy, which is a fixed plan chosen at the initial time and enforced irrevocably across all future periods, regardless of new information or changing market conditions.  On the other hand, an equilibrium strategy or time-consistent strategy is a plan that remains optimal at every period when re-evaluated from that period's perspective, which aligns with the concept of a subgame perfect equilibrium in dynamic games, where strategies are optimal in all "subgames" (i.e., at all points in time).

In this paper, we  study the equilibrium strategy, especially the extended HJB equation related to it.  This means that we need  to formulate the problem as a game and look for equilibrium.   It is indeed a fundamental result in game theory that equilibrium generally exists for mixed strategies rather than pure strategies in a broad class of games. Hence, it is better to consider  some notion of mixed strategy, namely, at each time, the players fix a probability of stopping and decide whether or not to stop according to this probability.  Such  mixed strategy is also called randomized stopping time in some literature, see, for example, Bayer et al. \cite{bayer2021randomized} and Dong \cite{dong2024randomized}. Here, we  model it as a doubly stochastic Poisson process. More precisely, let $\Theta$  be a random variable, which is exponentially distributed with unit intensity and independent of Brownian motion. Given a non-negative $\mathbb F$-adapted process $\{ \pi_s\}_{t\le s\le T}$, a random time $\tau$ is defined as 
$$
\tau:= \inf \left\{ s\in [t,T]: \int_t^s \pi_u du \ge \Theta\right\}  {\wedge} T,
$$
where we adopt the convention that the infimum of an empty set is infinity.   It represents the time that the player chooses to stop. Literally speaking,   it means that, conditioning on having not stopped before, the probability that  the player stops between $t$ and $t+dt$ is $\pi_t dt$.  Under this formulation, instead of choosing a  stopping time, the player  chooses the  intensity process $\pi$ to  optimize the  MV objective function.

From the results of doubly stochastic Poisson processes (see Jeanblanc et al.~\cite{jeanblanc2009mathematical} for details), one can compute that, for any function $\varphi$,
$$
\mathbb E \left[ \varphi(X_\tau) |\mathcal F_T\right]=\int_t^T \varphi(X_s)\pi_s e^{-\int_t^s \pi_u du}ds+\varphi(X_T)e^{-\int_t^T \pi_udu}.
$$
Then, we have that 
$$
\mathbb E_{t,x}[\varphi(X_\tau)]=\mathbb E_{t,x}[\mathbb E[\varphi(X_\tau)|\mathcal F_T]]=\mathbb E_{t,x}\left[\int_t^T \varphi(X_s)\pi_s e^{-\int_t^s \pi_u du}ds+\varphi(X_T)e^{-\int_t^T \pi_udu}\right].
$$
Thus, the  MV criteria can be rewritten as 
\begin{equation*}
\begin{split}
&\mathbb E_{t,x}\left[f(X_\tau)\right]-\frac{\gamma}{2}\text{Var}_{t,x}\left[f(X_\tau)\right]\\
=&\mathbb E_{t,x}\left[\int_t^T (f-\frac{\gamma}{2}f^2)(X_s)\pi_s e^{-\int_t^s \pi_u du}ds+(f-\frac{\gamma}{2}f^2)(X_T)e^{-\int_t^T \pi_u du}\right]\\
&+\frac{\gamma}{2}\left(\mathbb E_{t,x}\left[\int_t^T f(X_s)\pi_s e^{-\int_t^s \pi_u du}ds+f(X_T)e^{-\int_t^T \pi_u du}\right]\right)^2.
\end{split}
\end{equation*}
However, we consider a regularized   MV  reward, which is defined as the following
\begin{eqnarray}
			J^{\lambda}(t,x;\pi)&:=&\mathbb E_{t,x}\left[\int_t^T\left\{  (f-\frac{\gamma}{2}f^2)(X_s)\pi_s+\lambda H(\pi_s)\right\}e^{-\int_t^s \pi_udu}ds+(f-\frac{\gamma}{2}f^2)(X_T)e^{-\int_t^T \pi_udu}\right] \nonumber\\
		&&{}+\frac{\gamma}{2}\left(\mathbb E_{t,x}\left[ \int_t^T f(X_s)\pi_s e^{-\int_t^s \pi_udu}ds+f(X_T)e^{-\int_t^T \pi_udu}\right] \right)^2 \label{J^lambda},
	\end{eqnarray}
where the function $H$ is given by   $H(\pi)=\pi-\pi\log\pi$ and $\lambda$ is a positive constant.
The regularization term prevents the intensity taking values $0$ and $\infty$, which refers to continue and stop deterministically.  This provides mathematically tractability for the problem, which one can derive the related HJB equations.  The motivation comes from \cite{dong2024randomized} in the study of reinforcement learning method for optimal stopping problem, where $H$ is referred as unnormalized entropy to encourage randomness in the strategy.  

We  will let $\lambda$ tend to zero to go back to the original problem. Dong \cite{dong2024randomized} proves that the HJB equations  converge for optimal stopping problem. To our best knowledge, it is still an open problem for MV stopping problem, but, formally, we will see that the extended HJB equation converge to some  system that gives an equilibrium stopping time  for MV problem.

Since we focus on extended HJB equation, we shall restrict to the Markovian strategies, i.e., $\pi_s=\pi(s,X_s)$ for some deterministic function $\pi$.  Similar to \cite{bjork2017time}, one can define an equilibrium in the following sense. 
\begin{definition}
A strategy $\pi^*$ is called an equilibrium strategy, if for any $\varepsilon,v >0$ and $ t\in [0,T)$, define the perturbed policy $\pi^{\varepsilon,v}$  as 
\begin{equation*}
\pi^{\varepsilon,v}(s,x)=
\left\{
\begin{split}
&v,\text{ if $t \le s \le t+\varepsilon$};\\
&\pi^*(s,x)\text{ if $s>t+\varepsilon$},
\end{split}
\right.
\end{equation*}
and it holds that 
$$
\liminf_{\varepsilon \rightarrow 0} \frac{J^{\lambda}(t,x;\pi^*)-J^{\lambda}(t,x;\pi^{\varepsilon,v})}{\varepsilon} \ge 0, a.s.,
$$
for any initial state $x$.
\end{definition}

\section{Extended HJB Equation for Regularized Problem}
In this section, we derive the extended HJB equation for the regularized  MV problem. To this end, we first define an operator $\mathcal L$ as, for any smooth function $\varphi$, 
$$
(\mathcal L \varphi)(t,x):=\frac{1}{2}\mathop{\text{tr}}( \sigma\sigma^T(t,x)\partial_{xx}\varphi(t,x))+ b(t,x) \partial_x \varphi(t,x).
$$
We have the following verification theorem. 
\begin{theorem}\label{thm_equil} 
For a Markovian strategy $\pi^*=\pi^*(t,x)$,  let $(V^\lambda,g^\lambda)$ be a classical solution of the following parabolic system 
\begin{equation}\label{extended-HJB}
	\left\{
	\begin{split}
		&\partial_t V^\lambda+\mathcal LV^\lambda+\lambda\exp(-\frac{V^\lambda+\frac{\gamma}{2}(f-g^\lambda)^2-f}{\lambda})-\gamma |\sigma \partial_x g^\lambda|^2=0,\ V^\lambda(T,x)=f(x),\\
		& \partial_t g^\lambda+\mathcal L g^\lambda-\exp(-\frac{V^\lambda+\frac{\gamma}{2}(f-g^\lambda)^2-f}{\lambda})(g^\lambda-f)=0,\ g^\lambda(T,x)=f(x).
	\end{split}
	\right.
\end{equation}
Assume that $V^\lambda$,  $g^\lambda$, their derivatives (up to first order in $t$ and  second order in $x$) and $\pi$ are all continuous with polynomial growth in $x$, uniformly  in $t$. Then $\pi^*$ is an equilibrium strategy if and only if
\begin{equation}\label{condi_optim}
	\pi^*(t,x)=\exp(-\frac{V^\lambda+\frac{\gamma}{2}(f-g^\lambda)^2-f}{\lambda}).
\end{equation}
\end{theorem}
\begin{proof}
	The proof consists of two steps:
	\begin{enumerate}
		\item We start by showing that $g^\lambda(t,x)= \mathbb E_{t,x}[f(X_\tau)]$ and $V^\lambda(t,x)=J^\lambda(t,x;\pi^*)$;
		\item In the second step, we prove that $\pi^*$ is an equilibrium if and only if  \eqref{condi_optim} holds. 
	\end{enumerate}
For the first step, applying It\^o formula to $g^\lambda(s,X_s)e^{-\int_t^s \pi^*(u,X_u)du}$, we have that 
\begin{equation*}
\begin{split}
d\left(g^\lambda(s,X_s)e^{-\int_t^s \pi^*(u,X_u)du}\right)=&e^{-\int_t^s \pi^*(u,X_u)du}\left((\partial_t+\mathcal L)g^\lambda(s,X_s)-\pi^*(u,X_u)g^\lambda(u,X_u)\right)ds\\
&+\sigma \partial_x g^\lambda(s,X_s)e^{-\int_t^s \pi^*(u,X_u)du} dW_s.
\end{split}
\end{equation*}
From the growth assumption on the derivative of $g^\lambda$, it holds that the stochastic integral is a martingale. Thus, taking conditional expectation, one can verify that    
$$
g^\lambda(t,x)=\mathbb E_{t,x}\left[\int_t^T f(X_s)\pi^*(s,X_s)e^{-\int_t^s \pi^*(u,X_u) du}ds+f(X_T)e^{-\int_t^T \pi^*(u,X_u) du}\right],
$$
Then, define the function $h^\lambda$ as  $h^\lambda=V^\lambda-\frac\gamma2 (g^\lambda)^2$. It is easy to see that $h^\lambda$ solves
$$
(\partial_t +\mathcal L) h^\lambda(t,x)+\pi^*(t,x)(f-\frac{\gamma}{2}f^2-h^\lambda)(t,x)+\lambda H(\pi^*)=0,\ h^\lambda(T,x)=f(x)-\frac{\gamma}{2}f^2(x).
$$
Similarly, one can verify that 
$$
h^\lambda(t,x)=\mathbb E_{t,x}\left[\int_t^T \left\{(f-\frac{\gamma}{2}f^2)(X_s)\pi^*(s,X_s)+\lambda H(\pi_s) \right\}e^{-\int_t^s \pi^*(u,X_u) du}ds+f(X_T)e^{-\int_t^T\pi^*(u,X_u) du}\right].
$$
This implies that $J^\lambda(t,x;\pi^*)=h^\lambda(t,x)+\frac{\gamma}{2}g^\lambda(t,x)=V^\lambda(t,x)$.

Next, we prove the second step. For any $t$, $\varepsilon$ and $v$, consider the perturbed strategy $\pi^{\varepsilon,v}$. Since $\pi^{\varepsilon,v}$ coincides with $\pi^{*}$ after time $t+\varepsilon$, it holds that 
\begin{equation*}
\begin{split}
J^\lambda(t,x;\pi^{\varepsilon,v})=&\mathbb E_{t,x}\left[\int_t^{t+\varepsilon} \left\{(f-\frac{\gamma}{2}f^2)(X_s)v+\lambda H(v)\right\}e^{-v(s-t)} ds+h^\lambda(t+\varepsilon,X_{t+\varepsilon})e^{-v \varepsilon}\right]\\
&+\frac{\gamma}{2}\left(\mathbb E_{t,x}\left[\int_t^{t+\varepsilon}f(X_s)ve^{-v(s-t)}ds +g^\lambda(t+\varepsilon,X_{t+\varepsilon})e^{-v\varepsilon}\right]\right)^2.
\end{split}
\end{equation*}
Applying It\^o formula to $g^\lambda(s,X_{s})e^{-v(s-t)}$ and taking conditional expectation, we get that 
\begin{equation*}
\begin{split}
\mathbb E_{t,x}[g^\lambda(t+\varepsilon,X_{t+\varepsilon})e^{-v \varepsilon}]=g^\lambda(t,x)+\mathbb E_{t,x}\left[\int_t^{t+\varepsilon} (\partial_t+\mathcal L)g^\lambda(s,X_s)e^{-v(s-t)}-vg^\lambda(s,X_s)e^{-v(s-t)}ds \right].
\end{split}
\end{equation*}
Then, it holds that
\begin{equation*}
\begin{split}
&\mathbb E_{t,x}\left[\int_t^{t+\varepsilon}f(X_s)ve^{-v(s-t)}ds +g^\lambda(t+\varepsilon,X_{t+\varepsilon})e^{-v\varepsilon}\right]\\
=&g^\lambda(t,x)+\mathbb E_{t,x}\left[\int_t^{t+\varepsilon}(\partial_t+\mathcal L) g^\lambda(s,X_s)e^{-v(s-t)}+v(f-g^\lambda)(s,X_s)e^{-v(s-t)} ds\right]\\
=&g^\lambda(t,x)+\bigg((\partial_t+\mathcal L) g^\lambda(t,x)+v(f-g^\lambda)(t,x) \bigg)\varepsilon+o(\varepsilon),
\end{split}
\end{equation*}
where the last equality is  due to the fact that the related functions are continuous and $X$ also has continuous trajectories. Thus, we get that 
\begin{equation*}
\begin{split}
  &\left(\mathbb E_{t,x}\left[\int_t^{t+\varepsilon}f(X_s)ve^{-v(s-t)}ds +g^\lambda(t+\varepsilon,X_{t+\varepsilon})e^{-v\varepsilon}\right]\right)^2\\
  =&(g^\lambda)^2(t,x)+2 g^\lambda(t,x)\bigg((\partial_t+\mathcal L) g^\lambda(t,x)+v(f-g^\lambda)(t,x)\bigg)\varepsilon+o(\varepsilon).
  \end{split}
\end{equation*}
With a similar argument, one also have that 
\begin{equation*}
\begin{split}
&\mathbb E_{t,x}\left[\int_t^{t+\varepsilon} \left\{(f-\frac{\gamma}{2}f^2)(X_s)v+\lambda H(v)\right\}e^{-v(s-t)} ds+h^\lambda(t+\varepsilon,X_{t+\varepsilon})e^{-v \varepsilon}\right]\\
=&h^\lambda(t,x)+\bigg((\partial_t +\mathcal L) h^\lambda(t,x)+v(f-\frac{\gamma}{2}f^2-h^\lambda)(t,x)+\lambda H(v)\bigg)\varepsilon +o(\varepsilon).
\end{split}
\end{equation*}
Combining these estimations together, we get that
\begin{equation*}
 \begin{split}
J^\lambda(t,x;\pi^{\varepsilon,v})=&h^\lambda(t,x)+\frac\gamma2(g^\lambda)^2(t,x)+\bigg(
(\partial_t+\mathcal L)h^\lambda(t,x)+\gamma g^\lambda(t,x)(\partial_t+\mathcal L)g^\lambda(t,x)\\
&+v(f-\frac\gamma2 f^2-h^\lambda+\gamma g^\lambda f-\gamma (g^\lambda)^2)(t,x)+\lambda H(v) 
\bigg)\varepsilon+o(\varepsilon).
 \end{split}   
\end{equation*}
From the equations satisfied by $g^\lambda$ and $h^\lambda$, also noting that $J^\lambda(t,x;\pi^*)=h^\lambda(t,x)+\frac\gamma2(g^\lambda)^2(t,x)$, we derive that
\begin{equation*}
\begin{split}
&\mathop{\text{lim}}_{\varepsilon \rightarrow 0} \frac{J^{\lambda}(t,x;\pi^*)-J^{\lambda}(t,x;\pi^{\varepsilon,v})}{\varepsilon}\\
=&(v-\pi^*(t,x))(f-\frac\gamma2 f^2-h^\lambda+\gamma g^\lambda f-\gamma (g^\lambda)^2)(t,x)+\lambda(H(v)-H(\pi^*(t,x))).
\end{split}
\end{equation*}
Hence, $\pi^*$ is an equilibrium strategy if and only if
$$
\pi^*(t,x) \in \text{argmax}_{v} v(f-\frac{\gamma}{2}f^2-h^\lambda-\gamma (g^\lambda)^2+\gamma g^\lambda f)(t,x)+\lambda H(v),
$$
which implies that $\pi^*$ satisfies the optimality condition \eqref{condi_optim}.
\end{proof}

To our best knowledge, it is still a hard open problem to prove existence and/or uniqueness for solutions of an extended HJB system with a general assumption. In linear-quadratic mean-variance problem, one can reduce it to an ODE system and obtain a solution, see \cite{bjork2014mean}.   
For the solvability of \eqref{extended-HJB}, we give an existence result with additional technical assumptions on the coefficients and a small time interval.
\begin{theorem} \label{thm3.2}
 In addition to Assumption \ref{ass}, we further assume that the coefficients $b$,$\sigma$ and $f$ are uniformly bounded.  Then, for a sufficiently small $T$, \eqref{extended-HJB} admits a classical solution $(V^\lambda,g^\lambda)$.
\end{theorem}
\begin{proof}The solution is to be considered as a fixed point of a contraction mapping. For that purpose, 
denote by $\mathbb K$ the Banach space $C([0,T],C^1(\mathbb R^d))$ equipped with the norm $\|l\|_{\mathbb K}=\sup_{t,x}|l(t,x)|+\sup_{t,x}|\partial_x l(t,x)|$. For any $l \in \mathbb K$, define a mapping $F$ as $k=F(l)$ is the solution of the following system
\begin{equation}\label{decouple-HJB}
\left\{
\begin{split}
&\partial_t v+\mathcal Lv+\lambda\exp(-\frac{v+\frac{\gamma}{2}(f-l)^2-f}{\lambda})-\gamma |\sigma \partial_x l|^2=0, v(T,x)=f,\\
& \partial_tk+\mathcal L k-\exp(-\frac{v+\frac{\gamma}{2}(f-l)^2-f}{\lambda})(k-f)=0,k(T,x)=f.
\end{split}
\right.
\end{equation}
Let $B_m(0)$ be the ball in $\mathbb K$ centered at $0$ with radius $m=\|f\|_{\infty}+\|\partial_x f\|_{\infty}+1$. We are to show that, for a sufficiently small $T$, $F$ is a contraction from $B_m(0)$ into itself and, thus, admits a fixed point, which would  be the solution of \eqref{extended-HJB}. 

The proof consists of several steps:

{\it Step $1$. Well-posedness of first equation.} Choose any $N>0$, let $\xi_N$ be a smooth cutoff function such that $\xi_N(x)=x$, for $x\le N$; and $\xi_N(x)=N+1$, for $x\ge N+1$. Consider the following equation 
$$
\partial_t v+\mathcal Lv+\lambda\exp(\xi_N(-\frac{v+\frac{\gamma}{2}(f-l)^2-f}{\lambda}))-\gamma |\sigma \partial_x l|^2=0, v(T,x)=f.
$$
Noting that the third term is a bounded Lipschitz function of $v$, it admits a solution $v$. Lemma \ref{append_lemma} yields that 
$$
-C( \|f\|_{\infty}+\lambda \exp(\frac{N+1}{\lambda}))\le v\le C(Tm^2+\|f\|_{\infty}).
$$
Let us give a refined lower bound estimation independent of $N$, which implies that $v$ solves the first equation in \eqref{decouple-HJB} when $N$ is sufficiently large.
Denote by $\psi(x)=\sqrt{1+|x|^2}$.  One can compute that 
$$
D_x \psi=\frac{x}{\sqrt{1+|x|^2}}\text{ and }  D^2_x \psi= \frac{1}{\sqrt{1+|x|^2}}I -\frac{1}{\left(1+|x|^2\right)^{\frac12}} x\otimes x,
$$
and, thus, $\mathcal L \psi \le C$. Since $v$ is a bounded function, it holds that, for any $\varepsilon>0$, $v+\varepsilon\psi$ attains a minimum at some point $(t^*,x^*)$. If $t^*=T$, then the terminal condition implies that 
$$v\ge -\|f\|_\infty-\varepsilon \psi.
$$
If $t^*<T$, it holds that, at $(t^*,x^*)$,
$$
D^2_x v \ge -\varepsilon D^2_x \psi, D_x v=-\varepsilon D_x \psi,  \text{ and } \partial_t v \ge 0. 
$$
This yields that  $(\partial_t +\mathcal L)v(t^*,x^*) \ge-\varepsilon \mathcal L\psi \ge -C\varepsilon$. On the other hand, we have  
$$
(\partial_t +\mathcal L)v=\gamma |\sigma D_x l|^2-\lambda\exp(\xi_N(-\frac{v+\frac{\gamma}{2}(f-l)^2-f}{\lambda})).
$$
Combining these two inequalities, we get that 
$$
\lambda\exp(\xi_N(-\frac{v+\frac{\gamma}{2}(f-l)^2-f}{\lambda})) \le \gamma |\sigma D_x l|^2 +C\varepsilon,
$$
or, equivalently,
$$
\xi_N(\frac{f-v-\frac{\gamma}{2}(f-l)^2}{\lambda})\ \le \log \frac{\gamma|\sigma \partial_x l|^2+C\varepsilon}{\lambda}.
$$
This implies that,  at $(r^*,x^*)$,
$$
v \ge f-\lambda  \log \frac{\gamma|\sigma \partial_x l|^2+C\varepsilon}{\lambda}.
$$
Thus, we have a lower bound estimation
$$
v\ge -\|f\|_\infty +\lambda \log \frac{\gamma Cm^2+C\varepsilon}{\lambda}-\varepsilon \psi.
$$
Letting $\varepsilon$ go to $0$, we finally get that 
$$
v\ge -\|f\|_\infty +\lambda \log \frac{\gamma Cm^2}{\lambda}.
$$

{\it Step 2. Bound and derivative estimations for $k$.} Note that the second estimation of \eqref{decouple-HJB} is just a linear equation of $k$. The well-posedness is straight-forward and one obtains the following bound estimation from Lemma \ref{append_lemma} 
$$
\|k\|_{\infty} \le \|f\|_{\infty}.
$$
To give estimation of the derivative, note that, from lower bound estimation for $V$ in the first step, it holds that 
$$
\exp(-\frac{v+\frac\gamma2(f-l)^2-f}{\lambda}) \le (\frac{\gamma Cm^2}{\lambda})^{\frac{C}{\lambda}(1+\|f\|^2_\infty+m^2)}.
$$
Hence, using Lemma \ref{append_lemma} again, we have that 
$$
\|\partial_x k\|_{\infty} \le (1+\sqrt{T})\|\partial_x f\|_{\infty}+C(\sqrt{T}+T)(\frac{\gamma Cm^2}{\lambda})^{\frac{C}{\lambda}(1+\|f\|^2_\infty+m^2)}. 
$$
Then, for a sufficiently small $T$, $F$ is mapping from $B_m(0)$ into itself. 

{\it Step 3. Contraction of the mapping $F$.}
For any $h_i \in B_m(0)$ with $i=1,2$, let $(v_i,k_i)$ be the related solutions of \eqref{decouple-HJB}.  Denote by $\delta v=v_1-v_2$, $\delta k=k_1-k_2$, and $\delta l=l_1-l_2$. Then, we that $(\delta v,\delta k)$ satisfy 
\begin{equation}\label{delta-HJB}
\left\{
\begin{split}
&(\partial_t+\mathcal L)\delta v- \delta e(\delta v-\frac\gamma 2(2f-l_1-l_2)\delta l)-\gamma(|\sigma \partial_x l_1|^2-|\sigma \partial_x l_2|^2)=0, \delta v(T,x)=0,\\
& (\partial_t+\mathcal L) \delta k-\exp(-\frac{v_1+\frac{\gamma}{2}(f-l_1)^2-f}{\lambda})\delta k+(k_2-f)\delta e\frac{\delta v-\frac\gamma 2(2f-l_1-l_2)\delta l}{\lambda} =0,\delta k(T,x)=0.
\end{split}
\right.
\end{equation}
with 
$$
\delta e=\int_0^1 \exp(-\frac{v_1+\frac{\gamma}{2}(f-l_1)^2-f+s(\delta v-\frac\gamma 2(2f-l_1-l_2)\delta l)}{\lambda})ds.
$$
From the estimation in previous step, we know that  the related functions $v_i,k_i$ and $h_i$ are uniformly bounded. Thus, $\delta e$ is a bounded function. Moreover, one gets that 
$$
\|\delta e \frac \gamma2(2f-l_1-l_2)\delta l\|_\infty \le \frac\gamma 2\| \delta e(2f-l_1-l_2)\|_{\infty}\|\delta l\|_\infty \le C\|\delta l\|_\infty,
$$
and 
$$
\||\sigma D_x l_1|^2-|\sigma D_x l_2|^2 \|_\infty \le C\| \sigma D_x l_1+\sigma D_x l_2\|_\infty\| \sigma D_x l_1-\sigma D_x l_2\|_\infty \le C\| \sigma D_x l_1-\sigma D_x l_2\|_\infty.
$$
Then, using Lemma \ref{append_lemma} again, we have that there exists a constant $C$ depending on the coefficients, $m$ and $\lambda$ such that 
$$
\|\delta v\|_{\infty} \le CT(\|\delta l\|_{\infty}+\|\partial_x \delta l\|_{\infty}. 
$$
Furthermore, with a similar argument, 
$$
\|\delta k\|_\infty \le CT\|(k_2-f)\delta e\frac{\delta v-\frac\gamma 2(2f-l_1-l_2)\delta l}{\lambda}\|_\infty\le CT(\|\delta v\|_\infty+\|\delta l\|_\infty).,
$$
and
\begin{equation*}
\begin{split}
\|\partial_x \delta k\|_{\infty} \le& C(\sqrt{T}+T)(\|(k_2-f)\delta e\frac{\delta v-\frac\gamma 2(2f-l_1-l_2)\delta l}{\lambda}\|_\infty+\|\exp(-\frac{v_1+\frac{\gamma}{2}(f-l_1)^2-f}{\lambda})\delta k\|_\infty)\\
\le&C(\sqrt{T}+T)(\|\delta v\|_\infty+\|\delta l\|_\infty+\|\delta k\|_\infty).
\end{split}
\end{equation*}
Combining these estimations, we  see that $F$ is a contraction for a sufficiently small $T$.
\end{proof}

\section{Extended HJB Equation for Original Problem}
For the optimal stopping problem, i.e., $\gamma=0$,  \cite{dong2024randomized} proves that the first equation in \eqref{extended-HJB} converges to the variational inequality satisfied by the value function.  For  MV stopping problem, we would like follow the same procedure.  Unfortunately, at present, we can only formally deduce  the limiting equation. Below is a brief introduction. Assume that 
$(V^\lambda,g^\lambda)$ converge to some functions $(V,g)$ when $\lambda$ goes to $0$ Then, it is natural to conjecture  that  
$$\partial_t V+\mathcal LV-\gamma\left|\sigma\partial_x g \right|^2=\lim_{\lambda \rightarrow 0}\partial_t V^\lambda+\mathcal LV^\lambda-\gamma\left|\sigma\partial_x g^\lambda \right|^2\le 0.
$$
If $\partial_t V^\lambda+\mathcal LV^\lambda$ and $\partial_x g^\lambda$ are bounded uniformly in $\lambda$, then  so is $\lambda\exp(-\frac{V^\lambda+\frac{\gamma}{2}(f-g^\lambda)^2-f}{\lambda})$. This suggests that $V+\frac\gamma2(f-g)^2 \ge f$. Moreover, if $V+\frac\gamma2(f-g)^2 > f$, one has $\exp(-\frac{V^\lambda+\frac{\gamma}{2}(f-g^\lambda)^2-f}{\lambda})$ converges to $0$, which implies   $\partial_t V+\mathcal LV-\gamma\left|\sigma\partial_x g \right|^2=0$, and $\partial_t g+\mathcal Lg=0$. From \eqref{extended-HJB}, it is not clear what is satisfied by $g$ on the set $\{V+\frac\gamma2(f-g)^2=f\}$. But, from the probabilistic representation of $g^\lambda$, one also guess that $g(t,x)=\mathbb E_{t,x}[f(X_\tau)]$ for a stopping time $\tau$. If $\tau$ is the hitting time of the set $\{V+\frac\gamma2(f-g)^2=f\}$, Then, $g$ should equal to $f$ on that set.  In summary, \eqref{extended-HJB} formally converges to the following system
\begin{equation}\label{extended_VI}
\left\{
\begin{split}
&\min\left \{-\left(\partial_t V+\mathcal LV-\gamma\left|\sigma\partial_x g \right|^2\right),V+\frac{\gamma}{2}(f-g)^2-f\right\}=0,V(T,x)=f,\\
& \partial_t g+\mathcal Lg=0, \text{ on $\{V+\frac{\gamma}{2}(f-g)^2 >f\}$},\\
&g=f \text{ on $\{V+\frac{\gamma}{2}(f-g)^2 =f\}$}, g(T,x)=f(x).
\end{split}
\right.
\end{equation}

Now, let us assume that the above system admits a pair of solution $(V,g)$. The means that $(V,g)$ is a pair of continuous functions, second-order continuous differentiable in the region $\{V+\frac{\gamma}{2}(f-g)^2>f\}$ and satisfies \eqref{extended_VI}. Define  the set 
$$
\mathcal C=\{(t,x)| V+\frac{\gamma}{2}(f-g)^2>f\}.
$$
Construct the stopping time $\tau_{\mathcal C}$ as 
$$
\tau_{\mathcal C}=\inf\{s\ge t| (s,X_s) \notin \mathcal C\}.
$$
One can verify that 
$$
V(t,x)=\mathbb E_{t,x}\left[f(X_{\tau_{\mathcal C}})\right]-\frac\gamma2\text{Var}_{t,x}\left[f(X_{\tau_{\mathcal C}})\right], \text{ and } g(t,x)=\mathbb E_{t,x}\left[f(X_{\tau_{\mathcal C}})\right].
$$
Moreover, we also define 
$$
h(t,x):=V(t,x)-\frac{\gamma}{2}g^2(t,x)=\mathbb E_{t,x}\left[(f-\frac{\gamma}{2}f^2)(X_{\tau_{\mathcal C}})\right].$$

Next, we prove that these functions characterize a stopping policy that is an equilibrium in some sense. Before that, let us first introduce some definitions. For any $(t,x)\in \mathbb R \times\mathbb R^n$ and $r>0$, the parabolic cylinder $Q(t,x;r)$ is defined as\footnote{The definition of parabolic cylinders is different from that in text books of parabolic PDEs, see \cite{lieberman1996second}. The reason is that we consider PDEs with terminal conditions instead of initial conditions.} 
$$
Q(t,x;r):=\{(s,y)\in \mathbb R \times\mathbb R^n| \max\{|x-y|,(s-t)^{\frac12}\} \le r,s\ge t\}.
$$
For any set $\Omega$, the parabolic boundary $\mathcal P \Omega$ is defined as the set of all points $(t,x) \in \bar \Omega$ such that for any $\varepsilon>0$, $Q(t,x;\varepsilon)$ contains points not in $\Omega$. Finally, note that, since $\mathcal C$ is a open set, for any $(t,x) \in \mathcal C$, there exists $\varepsilon>0$ such that $Q(t,x;\varepsilon) \in \mathcal C$. 

The notion of equilibrium strategy is similar to that used in \cite{christensen2020time}. 
For any $v \ge 0$, let $N^{v}$ be a Poisson point process independent of the Brownian motion $W$. Its first jump time after $t$ is denoted as $\tau^v$, i.e.,
$$
\tau^v=\inf\left\{s \ge t| N_s-N_t=1\right\}.
$$
For any $\varepsilon$, define two stopping time 
$$
\tau^{\varepsilon}=\inf\left\{s\ge t| |X_s-X_t|\ge \varepsilon\right\} \wedge (t+\varepsilon) \wedge T,
$$
and 
$$
\tilde \tau_{\mathcal C}=\inf\left\{s \ge \tau^\varepsilon| (s,X_s) \notin \mathcal C\right\}.
$$
The   perturbation $\tau^{\varepsilon,v}_\mathcal C$ of $\tau_\mathcal C$  is defined as,
$$
\tau_{\mathcal C}^{\varepsilon,v}=1_{\{\tau^v \le \tau^\varepsilon \}} \tau^v +1_{\{\tau^v > \tau^\varepsilon \}}\tilde \tau_{\mathcal C}.
$$
\begin{theorem}\label{thm_stopping}
Assume that there exists a  solution $(V,g)$ of \eqref{extended_VI} such that the functions and their derivatives are polynomial-growth w.r.t. $x$ uniformly in $t$. For any  $(t,x) \in \mathcal C \cup (\mathcal C^c/\mathcal P \mathcal C^c)$, it holds that 
$$ \liminf_{\varepsilon\to0}\frac{J(t,x;\tau_{\mathcal C})-J(t,x;\tau_{\mathcal C}^{\varepsilon,v})}{\mathbb E_{t,x}[\tau^{\varepsilon}-t]}\ge 0.
$$
\end{theorem}
\begin{proof}
From the definition of $\tau_{\mathcal C}^{\varepsilon,v}$, it holds that  
\begin{equation*}
\begin{split}
&\mathbb E_{t,x}\left[f(X_{\tau_{\mathcal C}^{\varepsilon,v}}) \right]\\
=&\mathbb E_{t,x}\left[\int_t^{\tau^{\varepsilon}} f(X_s)ve^{-v(s-t)}ds+f(X_{\tilde \tau_{\mathcal C}}) e^{-v(\tau^{\varepsilon}-t)}\right]\\
=&\mathbb E_{t,x}\left[\int_t^{\tau^{\varepsilon}} f(X_s)ve^{-v(s-t)}ds+g(\tau^{\varepsilon},X_{\tau^{\varepsilon}}) e^{-v(\tau^{\varepsilon}-t)}\right].
\end{split}
\end{equation*}
Note that, for sufficiently small $\varepsilon$, $X_s$  stays in $\mathcal C$ or $\mathcal C^c$ for $s\in[t,\tau^\varepsilon]$ as we assume that $(t,x) \in \mathcal C \cup (\mathcal C^c/\mathcal P \mathcal C^c)$. Thus, one can apply It\^o formula to get that 
\begin{equation*}
\begin{split}
&\mathbb E_{t,x}\left[g(\tau^{\varepsilon},X_{\tau^{\varepsilon}}) e^{-v\tau^{\varepsilon}} \right]\\
=&g(t,x)+\mathbb E\left[ \int_t^{\tau^{\varepsilon}}(\partial_t+\mathcal L)g(s,X_s)e^{-v(s-t)}-vg(s,X_s)e^{-v(s-t)}ds\right].
\end{split}
\end{equation*}
Then, we have 
\begin{equation*}
\begin{split}
&\mathbb E_{t,x}\left[f(X_{\tau^{\varepsilon,v}_{\mathcal C}}) \right]\\
=&g(t,x)+\left((\partial_t+\mathcal L)g(t,x)+v(f-g)(t,x)\right)\mathbb E\left[\tau^{\varepsilon}-t\right]+o(\mathbb E\left[\tau^{\varepsilon}-t\right]),
\end{split}
\end{equation*}
which implies that 
\begin{equation*}
\begin{split}
&\left(\mathbb E_{t,x}\left[f(X_{\tau^{\varepsilon,v}_{\mathcal C}}) \right]\right)^2\\
=&g^2(t,x)+2g(t,x)\left((\partial_t+\mathcal L)g(t,x)+v(f-g)(t,x)\right)\mathbb E\left[\tau^{\varepsilon}-t\right]+o(\mathbb E\left[\tau^{\varepsilon}-t\right]).
\end{split}
\end{equation*}
Similarly, it holds that 
\begin{equation*}
\begin{split}
&\mathbb E_{t,x}\left[(f-\frac{\gamma}{2}f^2)(X_{\tau^{\varepsilon,v}_{\mathcal C}}) \right]\\
=&h(t,x)+\left((\partial_t+\mathcal L)(h)+v(f-\frac{\gamma}{2}f^2-h)\right)\mathbb E\left[\tau^{\varepsilon}-t\right]+o(\mathbb E\left[\tau^{\varepsilon}-t\right]).
\end{split}
\end{equation*}
Hence, recalling that $h=V-\frac\gamma2 g^2$,
\begin{equation*}
\begin{split}
&J(t,x;\tau_{\mathcal C}^{\varepsilon,v})-J(t,x;\tau_{\mathcal C})\\
=&((\partial_t +\mathcal L)h+\gamma g(\partial_t +\mathcal L)g+v(f-(V+\frac{\gamma}{2}(f-g)^2)))\mathbb E[\tau^{\varepsilon}-t]+o(\mathbb E[\tau^{\varepsilon}-t])\\
=&((\partial_t +\mathcal L)V-\gamma|\sigma \partial_x g|^2+v(f-(V+\frac{\gamma}{2}(f-g)^2)))\mathbb E[\tau^{\varepsilon}-t]+o(\mathbb E[\tau^{\varepsilon}-t]).
\end{split}
\end{equation*}
The first equation in \eqref{extended_VI} implies the desired result.
\end{proof}
For the case $(t,x) \in \mathcal P\mathcal C^c\backslash\{T\}\times\mathbb R$,  it turns out to be a very subtle problem. Thus, we focus on one dimensional case, i.e., $x\in \mathbb R$ and assume that the free boundary $\mathcal P\mathcal C^c\backslash\{T\}\times\mathbb R$ is locally Lipschitz continuous  with respect to time $t$. More precisely, there exists a small ball $Q(t,x;r)$ and a Lipschitz continuous curve $c$  and $r>0$ such that  $\mathcal C \bigcap Q(t,x;r)=\{(s,y)|t \le s \le t+r^{\frac12},|y-x|\le r, y\ge c(s)\}$.
\begin{theorem} \label{thm_stopping1}
In addition to the assumption in Theorem \ref{thm_stopping} and locally Lipschitz assumption on the free boundary, we further assume 
that $V$ is $C^1$ in $Q(t,x;r)$. If it holds that, for $(t,x) \in\mathcal P\mathcal C^c\backslash\{T\}\times\mathbb R$, 
\begin{equation}\label{ineq_boundary}
(\partial_t+\mathcal L )V(t,x+)+(\partial_t+\mathcal L )V(t,x-)\le \gamma \sigma^2(t,x)\left(\frac{\partial_x g(t,x+)+\partial_x g(t,x-)}{2}\right)^2,
\end{equation}
then we  have  
$$
\liminf_{\varepsilon \rightarrow 0}\frac{J(t,x;\tau_{\mathcal C})-J(t,x;\tau_{\mathcal C}^{\varepsilon,v})}{\mathbb E_{t,x}[\tau^{\varepsilon}-t]}\ge 0.
$$
\end{theorem}
\begin{proof}
From It\^o-Tanaka formula (see \cite{peskir2006optimal}), one can get that 
\begin{equation*}
\begin{split}
&g(\tau^{\varepsilon},X_{\tau^{\varepsilon}})e^{-(\tau^{\varepsilon}-t)}=g(t,x)+\int_t^{\tau^{\varepsilon}}\frac12\bigg( (\partial_t+\mathcal L)g(s,X_s+)e^{-v(s-t)}-vg(s,X_s+)e^{-v(s-t)}\\
&+(\partial_t+\mathcal L)g(s,X_s-)e^{-v(s-t)}-vg(s,X_s-)e^{-v(s-t)}\bigg)1_{\{X_s \neq c(s)\}}ds\\
&+\int_t^{\tau^{\varepsilon}}\frac12 \sigma( \partial_x g(s,X_s+)+ \partial_x g(s,X_s-))e^{-v(t-s)}dW_s\\
&+\frac{1}{2} \int_t^{\tau^{\varepsilon}}(\partial_x g(s,X_s+)-\partial_x g(s,X_s-))e^{-v(s-t)}1_{\{X_s=c(s)\}}dl^c_s, 
\end{split}
\end{equation*}
where $l^c_s$ is the local time of $X$ at the curve $c$. Then, we have that  
\begin{equation*}
\begin{split}
&\mathbb E_{t,x}\bigg[\int_t^{\tau^{\varepsilon}}\frac12\bigg( (\partial_t+\mathcal L)g(s,X_s+)e^{-v(s-t)}-vg(s,X_s+)e^{-v(s-t)}\\
&+(\partial_t+\mathcal L)g(s,X_s-)e^{-v(s-t)}-vg(s,X_s-)e^{-v(s-t)}\bigg)1_{\{X_s \neq c(s)\}}ds\bigg]\\
=&\frac{1}{2}\left(((\partial_t+\mathcal L)g(t,x+)+(\partial_t+\mathcal L)g(t,x-)-2vg(t,x)\right)\mathbb E[\tau^{\varepsilon}-t]+o(\mathbb E[\tau^{\varepsilon}-t]).
\end{split}
\end{equation*}
It is proved in Appendix \ref{append_esti} that
\begin{equation}\label{esti_square_of_local_time}
\left(\mathbb E_{t,x}[l^c_{\tau^{\varepsilon}}-l_t^c]\right)^2=\sigma^2(t,x)\mathbb E_{t,x}[\tau^{\varepsilon}-t]+o(\mathbb E_{t,x}[\tau^{\varepsilon}-t]),
\end{equation}
which further  implies that 
\begin{equation}
\begin{split}
&\left( \mathbb E_{t,x}[\int_t^{\tau^{\varepsilon}}(\partial_x g(s,X_s+)-\partial_x g(s,X_s-))e^{-v(s-t)}1_{\{X_s=c(s)\}}dl^c_s] \right)^2\\
=&\left( \partial_x g(t,x+)-\partial_x g(t,x-)\right)^2\sigma^2(t,x)\mathbb E_{t,x}[\tau^{\varepsilon}-t]+o(\mathbb E_{t,x}[\tau^{\varepsilon}-t]).
\end{split}
\end{equation}
Hence,
\begin{equation*}
\begin{split}
  &(\mathbb E_{t,x}[e^{-v\tau^{\varepsilon}}g(\tau^{\varepsilon},X_{\tau^{\varepsilon}})] )^2=g^2(t,x)+\bigg( ((\partial_t+\mathcal L)g(t,x+)+(\partial_t+\mathcal L)g(t,x-))g(t,x)-2vg^2_1(t,x)\\
&+\left(\frac{\partial_x g(t,x+)-\partial_x g(t,x-)}{2}\right)^2\sigma^2(t,x)\bigg)\mathbb E_{t,x}[\tau^{\varepsilon}-t]\\
  &+g(t,x)\mathbb E_{t,x}[\int_t^{\tau^{\varepsilon}}(\partial_x g(s,X_s+)-\partial_x g(s,X_s-))e^{-v(s-t)}1_{\{X_s=c(s)\}}dl^c_s] +o(\mathbb E_{t,x}[\tau^{\varepsilon}-t]).
  \end{split}
\end{equation*}
Similarly,
\begin{equation*}
\begin{split}
&\mathbb E_{t,x}[e^{-v\tau^{\varepsilon}}h(\tau^{\varepsilon},X_{\tau^{\varepsilon}})]\\
=&h(t,x)+\frac{1}{2}\left( (\partial_t+\mathcal L)h(t,x+)+(\partial_t+\mathcal L)h(t,x-)-vh(t,x)\right)\mathbb E_{t,x}[\tau^{\varepsilon}-t]\\
&+\frac{1}{2}\mathbb E_{t,x}[ \int_t^{\tau^{\varepsilon}}(\partial_x h(s,X_s+)-\partial_x h(s,X_s-))e^{-v(s-t)}1_{\{X_s=c(s)\}}dl^c_s].
\end{split}
\end{equation*}
Since $V=h+\frac{\gamma}{2}g^2$ is $C^1$ across $c$, it holds that 
\begin{equation*}
\begin{split}
&\mathbb E_{t,x}[ \int_t^{\tau^{\varepsilon}}(\partial_x h(s,X_s+)-\partial_x h(s,X_s-))e^{-v(s-t)}1_{\{X_s=c(s)\}}dl^c_s\\
&+\gamma g(t,x)\mathbb E_{t,x}[\int_t^{\tau^{\varepsilon}}(\partial_x g(s,X_s+)-\partial_x g(s,X_s-))e^{-v(s-t)}1_{\{X_s=c(s)\}}dl^c_s]=o(\mathbb E_{t,x}[\tau^{\varepsilon}-t]).
\end{split}
\end{equation*}
Then,
\begin{equation*}
\begin{split}
&\frac{J(t,x;\tau_{\mathcal C}^{\varepsilon,v})-J(t,x;\tau_{\mathcal C})}{\mathbb E_{t,x}[\tau^{\varepsilon}-t]}\\
=&\frac{1}{2}\left((\partial_t+\mathcal L)(V-\frac{\gamma}{2}g^2)(t,x+)+\gamma g(\partial_t+\mathcal L)g(t,x+)\right)\\
&+\frac{1}{2}\left((\partial_t+\mathcal L)(V-\frac{\gamma}{2}g^2)(t,x-)+\gamma g(\partial_t+\mathcal L)g(t,x-)\right)\\
&+\frac{\gamma}{2} \left(\frac{\partial_x g(t,x+)-\partial_x g(t,x-)}{2}\right)^2\sigma^2(t,x)+o(1)\\
=&\frac{1}{2}\bigg((\partial_t+\mathcal L)V(t,x+)-\gamma|\sigma \partial_x g|^2(t,x+)\\
&+(\partial_t+\mathcal L)V(t,x-)-\gamma|\sigma \partial_x g|^2(t,x-)\bigg)\\
&+\frac{\gamma}{2} \left(\frac{\partial_x g(t,x+)-\partial_x g(t,x-)}{2}\right)^2\sigma^2(t,x)+o(1),\\
\end{split}
\end{equation*}
which gives the desired result according  the assumption of the theorem.
\end{proof}
\begin{remark}
It is a well-known result in the theory of optimal stopping that the value function is $C^1$ across the free boundary. However, one can not expect that $g$ is also $C^1$, see the example given in Subsection \ref{subsec:1}. Note that a similar condition is also given in \cite{christensen2020time}. 
\end{remark}
\section{Further Discussions}
\subsection{Infinite Horizon Case}\label{subsec:1}
In this subsection, we consider an infinite horizon example that one can give an explicit solution.  Let $X$ be a geometric Brownian motion
$$
dX_t=\mu X_tdt+\sigma X_t dW_t.
$$
Consider the  MV stopping problem 
$$
J(x;\tau)=\mathbb E_x\left[X_\tau\right]-\frac{\gamma}{2}\text{Var}_x\left[X_\tau\right]. 
$$
For this infinite horizon case, one can have an elliptic system similar to \eqref{extended_VI}
\begin{equation}\label{extend_ellip}
\left\{
\begin{split}
&\min\left \{-\left(\mathcal L v-\frac{\gamma}{2}\left|\sigma x\partial_x g\right|^2\right),v+\frac{\gamma}{2}(f-g)^2-f\right\}=0,\\
& \mathcal Lg=0, \text{ on $\{v+\frac{\gamma}{2}(f-g)^2 >f\}$},\\
&g=f \text{ on $\{v+\frac{\gamma}{2}(f-g)^2 =f\}$},
\end{split}
\right.
\end{equation}
with $\mathcal L=\frac{1}{2} \sigma^2x^2 \frac{\partial^2}{\partial x^2}+\mu x\frac{\partial }{\partial x}$. Denote by $\rho=\frac{2\mu}{\sigma^2}$ and assume that $\rho \in (0,1/2)$. Let $b=\frac{2\rho}{\gamma (1-\rho)}$. Let us check that
\begin{equation*}
V(x)=\left\{
\begin{split}
&(1-\frac{\gamma}{2}b)b^{\rho}x^{1-\rho}+\frac{\gamma}{2} b^{2\rho}x^{2-2\rho},\text{ for $x<b$},\\
&x,\text{ for $x\ge b$},
\end{split}
\right.
\end{equation*}
and
\begin{equation*}
g(x)=\left\{
\begin{split}
&b^{\rho}x^{1-\rho},\text{ for $x<b$},\\
&x,\text{ for $x\ge b$},
\end{split}
\right.
\end{equation*}
satisfy \eqref{extend_ellip}.  For $x<b$, it holds that 
$$
\partial_x V=(1-\rho)(1-\frac\gamma 2 b)b^\rho x^{-\rho}+\gamma(1-\rho )b^{2\rho}x^{1-2\rho},
$$
$$
\partial_{xx} V=-\rho (1-\rho)(1-\frac\gamma 2 b)b^\rho x^{-\rho-1}+\gamma(1-\rho )(1-2\rho)b^{2\rho}x^{-2\rho},
$$
and $\partial_x g=(1-\rho b^\rho x^{-\rho})$. Thus, it is easy to check that, for $x<b$, $\mathcal L V-\frac{\gamma}{2}\left|\sigma x \partial_x g\right|^2=0$ and $\mathcal L g=0$.  When $x \ge b$, it is easy to see that $\partial_x V=\partial_x g=1$ and $\partial_{xx}V=0$. Hence,
$$
\mathcal L V-\frac{\gamma}{2}\left|\sigma x \partial_x g\right|^2=\mu x-\frac\gamma2\sigma^2 x^2=\frac\gamma2\sigma^2 x(\frac{\rho}{\gamma}-x) \le 0.
$$
Now let us verify that $V+\frac{\gamma}{2}(x-g)^2 >x$ for $x<b$. Direct computation yields that 
\begin{equation*}
\begin{split}
    &V+\frac{\gamma}{2}(x-g)^2\\
    =&(1-\frac{\gamma}{2}b)b^\rho x^{1-\rho}+\frac{\gamma}{2}b^{2\rho}x^{2-2\rho}+\frac{\gamma}{2}(x-b^\rho x^{1-\rho})^2\\
    =&\frac{\gamma}{2}x^2+\gamma b^{2\rho}x^{2-2\rho}-\gamma b^{\rho}x^{2-\rho}+(1-\frac{\gamma}{2}b)b^\rho x^{1-\rho}\\
    =&x\left(\frac{\gamma}{2} x+\gamma b^{2\rho}x^{1-2\rho} -\gamma b^{\rho}x^{1-\rho}+(1-\frac{\gamma}{2}b)b^{\rho}x^{-\rho}\right).
\end{split}
\end{equation*}
Thus, we have to show that $\frac{\gamma}{2} x+\gamma b^{2\rho}x^{1-2\rho} -\gamma b^{\rho}x^{1-\rho}+(1-\frac{\gamma}{2}b)b^{\rho}x^{-\rho}>1$ for $x<b$. For that purpose, define a function $\kappa(z)=\frac\gamma 2b z^{2\rho-1}+(1-\frac\gamma 2b)z^\rho$. Let us find its infimum on $[1,\infty)$. Taking derivative with respect to $z$, we see that $\kappa'(z)=z^{\rho-1}((2\rho-1)\frac\gamma 2 b z^{\rho-1}+(1-\frac \gamma 2 b)\rho)$. Since $\rho\le \frac12$ and $z\ge 1$, it holds that 
$$
(2\rho-1)\frac\gamma 2 b z^{\rho-1}+(1-\frac \gamma 2 b)\rho\ge (2\rho-1)\frac\gamma 2 b+ (1-\frac \gamma 2 b)\rho=0,
$$
where we use the fact that $b=\frac{2\rho}{\gamma(1-\rho)}$.
This implies that $\kappa$ is strictly increasing on $[1,\infty)$ and taking minimum at $z=1$, which equals to $1$. Then,  for $x<b$, 
\begin{equation*}
\begin{split}
&\frac{\gamma}{2} x+\gamma b^{2\rho}x^{1-2\rho} -\gamma b^{\rho}x^{1-\rho}+(1-\frac{\gamma}{2}b)b^{\rho}x^{-\rho}\\
=&\frac{\gamma}{2}x(1-b^\rho x^{-\rho})^2+\frac\gamma 2b^{2\rho}x^{1-2\rho}+(1-\frac\gamma2b)b^\rho x^{-\rho}\\
\ge &\kappa(bx^{-1})>1.
\end{split}
\end{equation*}
Thus, $(V,g)$ is a solution of  \eqref{extend_ellip}. On $x=b$, we can check that  
$$
\mathcal L V(b-)=\frac\gamma2 \sigma^2 b^2(1-\rho)^2 , \mathcal LV(b+)=\mu b, \partial_x g(b-)=1-\rho,\text{ and } \partial_xg(b+)=1. 
$$
Then,
$$
\mathcal L V(b-)+\mathcal LV(b+)=\frac\gamma2 \sigma^2 b^2(1-\rho)^2+\mu b=\sigma^2b(\frac \gamma 2 b(1-\rho)^2+\rho)=\sigma^2b\rho(2-\rho),
$$
and 
$$
\gamma \sigma^2b^2\left(\frac{1+1-\rho}{2}\right)^2=\sigma^2b(2-\rho)\frac{(2-\rho)\gamma b}{4}=\sigma^2b(2-\rho)\rho\frac{2-\rho}{2-2\rho}>\sigma^2b(2-\rho)\rho. 
$$
Thus,  the stopping time is an equilibrium.

Moreover, one can also compute that $V'(b+)=V'(b-)$, $g'_1(b-)=1-\rho$ and $g'_1(b+)=1$. This suggests that one can expect $V$ to be $C^1$, but $g$  not $C^1$, so  one can not assume that $g$ is $C^1$ across the free boundary.
\subsection{Discrete Time Approximation}
For many time-inconsistent problems, the equilibrium solution in continuous time is regarded as the limit of its counterpart in discrete-time models. This is the logical structure of the derivation for extended HJB equation, see \cite{bjork2021time}.  For that reason, we would like to check that whether one can get the same equation by considering the discrete time  MV stopping problem. 

Let $\Delta t=\frac{T}{N}$ and $t_k=k\Delta t$ for $k=0,1,...,N$. We assume that one can only stop at these time points. Recursively  define a sequence of stopping times $\{\tau_i\}_{i=0,1,...,N}$ as the follows. Set $\tau_N=t_N$, $U(N,x)=f(x)$, and $V(N,x)=f(x)$. For $i=N-1,N-2,...,0$, define 
$$
U(i,x):=\mathbb E_{t_i,x}[f(X_{\tau_{i+1}})]-\frac{\gamma}{2} \text{Var}_{t_i,x}[f(X_{\tau_{i+1}})],
$$
\begin{equation*}
\tau_i:=\left\{
\begin{matrix}
t_i,\text{ if $f(X_{t_i}) \ge U_{i}(X_{t_i})$,}\\
\tau_{i+1},\text{if $f(X_{t_i}) < U_{i}(X_{t_i})$,}
\end{matrix}
\right. 
\end{equation*}
and 
$$
V(i,x):=\mathbb E_{t_i,x}[f(X_{\tau_{i}})]-\frac{\gamma}{2} \text{Var}_{t_i,x}[f(X_{\tau_{i}})].
$$
Moreover, we also define $g(i,x):=\mathbb E_{t_i,x}[f(X_{\tau_{i}})]$. The motivation of these definitions are the following.  We view the  MV stopping problem from a game-theoretic perspective as  a non-cooperative game. We have one player at each time point  $t_n$, who can only choose the stopping decision at $t_n$. The stopping time $\tau_i$ represent the time the process being stopped conditioned on its has not been stopped before $t_i$. Player $n$ has two options: to stop immediately or to continue. The reward of the first option is $f(X_{t_i})$. Choosing to continue, the process  is stopped at time $\tau_{i+1}$ and the expected reward is given  by $U(i,X_{t_i})$. Thus, the strategy of player $n$ is to decide to stop  when $f(X_{t_i}) \ge U(i,X_{t_i})$.

Now let us check the equation satisfied by $V_i$. It is easy to see that $V(i,x) \ge f(x) $. If $V(i,x)>f(x)$, it implies that $\tau_i=\tau_{i+1}$. In this case, it holds that 
$$
V(i,x)=U(i,x)=\mathbb E_{t_i,x}[V(i+1,X_{t_{i+1}})]-\frac{\gamma}{2}\text{Var}_{t_i,x}[g(i+1,X_{t_{i+1}})].
$$
Combining two situations, we have 
$$
\min\left\{V(i,x)-\left( \mathbb E_{t_i,x}[V(i+1,X_{t_{i+1}})]-\frac{\gamma}{2}\text{Var}_{t_i,x}[g(i+1,X_{t_{i+1}})]\right),V(i,x)-f(x)\right\}=0.
$$
Now we let  $\Delta t$ go to zero. Formally, one see that 
\begin{equation*}
 \mathbb E_{t_i,x}[V(i+1,X_{t_{i+1}})]-V(i,x)=(\partial_t V+\mathcal L V)(t_i,x)\Delta t+o(\Delta t),
\end{equation*}
and 
\begin{equation*}
	\begin{split}
&\text{Var}_{t_i,x}[g(i+1,X_{t_{i+1}})]=\mathbb E_{t_i,x}[g^2(i+1,X_{t_{i+1}})]-\left(\mathbb E_{t_i,x}[g(i+1,X_{t_{i+1}})]\right)^2\\
=&g^2(i,x)+(\partial_t+\mathcal L)g^2\Delta t-\left(g(i,x)+(\partial_t+\mathcal L)g\Delta t \right)^2  o(\Delta t)\\
=&(\partial_t+\mathcal L)g^2\Delta t-g(\partial_t+\mathcal L)g\Delta t+o(\Delta t)\\
=&|\sigma \partial_x g|^2\Delta t+o(\Delta t).
\end{split}
\end{equation*}
Hence, we get the following system
\begin{equation}\label{extended_VI_classical}
\left\{
\begin{split}
&\min\left \{-\left(\partial_t V+\mathcal LV-\gamma\left|\sigma\partial_x g \right|^2\right),V-f\right\}=0,V(T,x)=f,\\
& \partial_t g+\mathcal Lg=0, \text{ on $\{V >f\}$},\\
&g=f \text{ on $\{V =f\}$}, g(T,x)=f(x).
\end{split}
\right.
\end{equation}
Note that this equation is different from \eqref{extended_VI} with the condition $V+\frac\gamma2(f-g)^2\ge f$ replaced by $V\ge f$. Clearly $V\ge f$ implies $V+\frac\gamma2(f-g)^2\ge f$. Moreover, on $\{V=f\}$, we also have $g=f$ and hence, $V+\frac\gamma2(f-g)^2=f$. Thus, the solution of the above system also satisfies 
\eqref{extended_VI}.  In this sense, it seems that \eqref{extended_VI} is a more general equation. However, it remains unclear whether there exists a solution that satisfies  \eqref{extended_VI}, but do not satisfy \eqref{extended_VI_classical}. 
\subsection{General Time-inconsistent Problems}
To explain the additional quadratic term in the condition $V+\frac\gamma2(f-g)^2\ge f$ of \eqref{extended_VI},  let us consider a more general time-inconsistent problem in which the player choose stopping time to maximize the following functional
$$
\mathbb E_{t,x}[f(X_\tau)]+G(\mathbb E_{t,x}[k(X_\tau)]). 
$$
We just give some formal arguments for illustration. First, we consider the regularized problem. For an equilibrium $\pi^*$, define $g(t,x)=\mathbb E_{t,x}[k(X_\tau)]$ and $h(t,x)=\mathbb E_{t,x}[f(X_\tau)+\int_0^\tau \lambda H(\pi^*_s)ds]$. Then, it holds that $(g,h)$ solves
\begin{equation*}
\left\{
\begin{split}
&(\partial_t+\mathcal L)h+\lambda H(\pi^*)+\pi^*(f- h)=0,g_2(T)=f,\\
&(\partial_t+\mathcal L)g+\pi^*(k-g)=0,g(T)=k.
\end{split}
\right. 
\end{equation*}
For a purterbed strategy $\pi^{\varepsilon,v}$, one can get that 
$$
\mathbb E_{t,x}[k(X_{\pi^{\varepsilon,v}})]=g(t,x)+\left((\partial_t+\mathcal L)g+v(k-g)\right)\varepsilon+o(\varepsilon),
$$
and
\begin{equation*}
	\begin{split}
		&\mathbb E_{t,x}[f(X_{\pi^{\varepsilon,v}})+\int_0^{\tau^{\varepsilon,v}} \lambda H(\pi^{\varepsilon,v}_s)ds]\\
		=&h(t,x)+\left((\partial_t+\mathcal L)h+v(f-h)+\lambda H(v)\right)\varepsilon+o(\varepsilon).
	\end{split}
\end{equation*}
Moreover, it holds that 
$$
G(\mathbb E_{t,x}[k(X_{\pi^{\varepsilon,v}})])=G(g)+G'(g)\left((\partial_t+\mathcal L)g+v(k-g)\right)\varepsilon+o(\varepsilon).
$$
Then, one can show that $\pi^*$ is an equilibrium if and only if
$$
\pi^*=\exp(-\frac{h+G'(g)g-f-G'(g)k}{\lambda}).
$$
Note that the value function $V=h+G(g)$ and satisfies 
$$
(\partial_t+\mathcal L)V+\lambda \exp(-\frac{h+G'(g)g-f-G'(g)k}{\lambda})+\mathcal H_G(g)=0,
$$
where the operator $\mathcal H_G$ is defined as 
$$
\mathcal H_G(\varphi)=G'(\varphi)\mathcal L\varphi-\mathcal  L G(\varphi). 
$$
We also see that 
$$
h+G'(g)g-f-G'(g)k=V-(f+G(k))+G(k)-G(g)-G'(g)(k-g).
$$
Denote by $\Delta_G(k,g)=G(k)-G(g)-G'(g)(k-g)$. Then, when letting $\lambda$ go to zero, $V$ should converge to the following variational inequality
$$
\min\left\{-\left((\partial_t+\mathcal L)V+H_G(g)\right), V-(f+G(k))+\Delta_G(k,g) \right\}=0.
$$
Note that the condition is written as  $V+\Delta_G(k,g)\ge (f+G(g))$.  For  MV case, i.e., $G(k)=\frac\gamma2 k^2$, $\Delta_G(k,g)=\frac\gamma2 (k-g)^2$, which is the same condition in \eqref{extended_VI}. Moreover, define the set 
$$
\mathcal C=\{(t,x)|V+\Delta_G(k,g)>f+G(k)\},
$$
and the stopping time $\tau_{\mathcal C}$ as 
$$
\tau_{\mathcal C}=\inf\{s>t|(s,X_s) \notin \mathcal C\}. 
$$
Then, one can  formally verify that $\tau_{\mathcal C}$ in the same sense as that in the statement of Theorem \ref{thm_stopping}.

From the previous discussion, we see that the appearance of  the additional term $\Delta_G(k,g)$ is due to two factors.  One is that  $G$ is a non-linear function, which makes the problem time-inconsistent. The other is that, when equilibrium is under consideration, the perturbation lies within the family of relaxed strategies rather than  pure strategies. This coincides with the fact in the game theory that pure strategy equilibrium is different from mixed strategy equilibrium.      

\section{{Conclusions}}

This paper systematically investigates the MV  optimal stopping problem--a time-inconsistent stochastic optimization problem-by developing a vanishing regularization method and deriving the corresponding extended HJB equation.  More precisely, to tackle the mathematical intractability of direct equilibrium analysis, we introduce a regularized problem, which enables rigorous derivation of the equilibrium strategy and the associated extended HJB equation. Then, letting $\lambda \rightarrow 0$, i.e. vanishing regularization, we formally recover a system of parabolic variational inequalities \eqref{extended_VI} for the original MV problem. This system characterizes equilibrium stopping times and includes a key quadratic term  $\frac{\gamma}{2}(f-g)^2$--a distinction from classical optimal stopping, where stopping conditions depend only on comparing the value function to the instantaneous reward. By extending the analysis to general time-inconsistent problems, we demonstrate that the additional term in our stopping condition arises from the non-linearity of the objective (responsible for time inconsistency) and the use of mixed strategies.

This work provides a rigorous mathematical foundation for MV stopping problems, with potential applications in financial decision-making (e.g., asset sale timing, portfolio liquidation) and statistical inference (e.g., risk-aware hypothesis testing stopping rules).  However, there are still some open problems left for future research, including
proving rigorous convergence of the regularized extended HJB system to the limiting variational inequality and developing numerical algorithms to compute equilibrium stopping rules for practical applications.

\bibliographystyle{siam}
\bibliography{references}
\appendix
\section{Proof of \eqref{esti_square_of_local_time}} \label{append_esti}
Define another stopping time $\tilde \tau^{\varepsilon}$ as 
$$
\tilde \tau^{\varepsilon}=\inf \{s\ge t||X_s-X_t|\ge \varepsilon\}.
$$
Note that $\tau^{\varepsilon} $ coincides to $\tilde \tau^{\varepsilon}$ on the set $\{\tilde \tau^{\varepsilon}<t+\varepsilon\}$. Then, we claim that 
\begin{equation}\label{tau_esti}
  \lim_{\varepsilon\rightarrow 0}\frac{\varepsilon^2}{\mathbb E_{t,x}[\tilde \tau^{\varepsilon}-t]}=\sigma^2(t,x).  
\end{equation}
To prove the claim, we follow the same argument as that in the proof of Lemma 5.5 in \cite{christensen2020time}. For any $a>\frac{1}{\sigma^2(t,x)}$, define a function $F$ as $F(t,y)=a(y-x)^2-t$. It is easy to see that 
$$
(\partial_t+\mathcal L)F(s,y)=2a\mu(s,y)(y-x)+a\sigma^2(s,y)-1,
$$
which is  greater than $0$ for $\{(s,y)||y-x|\le \varepsilon,s-t\le \varepsilon\}$ with a sufficiently small $\varepsilon$. Thus, applying It\^o formula to $F(s,X_s)$ and taking conditional expectation, it holds that 
$$\mathbb E_{t,x}[\tilde \tau^{\varepsilon}-t] \le a \mathbb E_{t,x}[|X_{\tilde \tau^{\varepsilon}}-x|^2]=a\varepsilon^2.$$
Similarly, one also has 
$$\mathbb E_{t,x}[\tilde \tau^{\varepsilon}-t] \ge a\varepsilon^2,$$
for $a<\frac{1}{\sigma^2(t,x)}$. Combining these two estimations  give the claim \eqref{tau_esti}.
Using Chebyshev's inequality, we also have that 
\begin{equation}\label{prob_esti}
    P(\tilde \tau^{\varepsilon}>t+\varepsilon) \le \frac{\mathbb E_{t,x}[\tilde \tau^{\varepsilon}-t]}{\varepsilon}=O(\varepsilon).
\end{equation}

The argument for \eqref{tau_esti} also yields that 
$$ \lim _{h\rightarrow 0}\frac{\mathbb E_{t,x}[|X_{\tau^{\varepsilon}}-x|^2]}{\mathbb E_{t,x}[\tau^{\varepsilon}-t]}=\sigma^2(t,x).
$$
Note that, using \eqref{prob_esti}, 
$$
\varepsilon^2\ge\mathbb E_{t,x}[|X_{\tau^{\varepsilon}}-x|^2]\ge \varepsilon^2P(\tau^{\varepsilon}<t+\varepsilon)=\varepsilon^2-\varepsilon^2P(\tau^{\varepsilon}=t+\varepsilon)=\varepsilon^2+o(\varepsilon^2).
$$
Hence, we get that 
\begin{equation}\label{est1}
 \mathbb E_{t,x}[\tau^{\varepsilon}-t]=\frac{1}{\sigma^2(t,x)}\varepsilon^2+o(\varepsilon).   
\end{equation}

Let $k(t,y):=|c(t)-y|$. From It\^o-Tanaka formula, one has that
\begin{equation*}
\begin{split}
k(\tau^{\varepsilon},X_{\tau^{\varepsilon}})=&\int_t^{\tau^{\varepsilon}}\frac{1}{2} ((\partial_t+\mathcal L)k(s,X_s+)+(\partial_t+\mathcal L)k(s,X_s-))ds\\
&+\int_t^{\tau^{\varepsilon}}\frac12\sigma(\partial_yk(s,X_s+)+\partial_yk(s,X_s-))dW_s\\
&+\int_t^{\tau^{\varepsilon}}\frac12(\partial_y k(s,X_s+)-\partial_yk(s,X_s-))1_{\{X_s=c(s)\}}dl_s^c. 
\end{split}
\end{equation*}
Note that, for $X_s=c(s)$, $\partial_y k(s,X_s+)=1$ and $\partial_y k(s,X_s+)=-1$. Thus, 
\begin{equation}\label{est2}
  (\mathbb E_{t,x}[k(\tau^{\varepsilon},X_{\tau^{\varepsilon}})])^2=(\mathbb E_{t,x}[l^c_{\tau^{\varepsilon}}-l^c_t])^2+o(\mathbb E_{t,x}[\tau^{\varepsilon}-t]).  
\end{equation}
When $t$ is fixed, for sufficienly small $\varepsilon$,
\begin{equation*}
\begin{split}
    \mathbb E_{t,x}[k(\tau^{\varepsilon},X_{\tau^{\varepsilon}})]=&\mathbb E_{t,x}[|x+\varepsilon-c(\tau^{\varepsilon})|1_{\{X_{\tau^{\varepsilon}=x+\varepsilon}\}}]+\mathbb E_{t,x}[|x-\varepsilon-c(\tau^{\varepsilon})|1_{\{X_{\tau^{\varepsilon}=x-\varepsilon}\}}]\\
    &+\mathbb E_{t,x}[|X_{\tau^{\varepsilon}}-c(\tau^{\varepsilon})|1_{\{\tau^{\varepsilon}=t+\varepsilon\}}].
    %=&p_{+}|x+\varepsilon-x|+p_{-}|x-\varepsilon-x|+\mathbb E[|X_{\tau^{\varepsilon}}-x|1_{\{\tau^{\varepsilon}=t+\varepsilon\}}]+O(\mathbb E[\tau^{\varepsilon}-t])\\
    %=&h+O(\mathbb E[\tau^{\varepsilon}-t])+O(h^2)
\end{split}
\end{equation*}
We estimate the right hand side term by term. Estimations for the first two term
are similar. It holds that 
\begin{equation*}
\begin{split}
\mathbb E_{t,x}[|x\pm \varepsilon-c(\tau^{\varepsilon})|1_{\{X_{\tau^{\varepsilon}=x\pm \varepsilon}\}}]&\le \mathbb E_{t,x}[\left(|x\pm \varepsilon-c(t)|+|c(t)-c(\tau^{\varepsilon})|\right)1_{\{X_{\tau^{\varepsilon}=x\pm \varepsilon}\}}]\\
&\le p_{\pm}\varepsilon+C\mathbb E_{t,x}[\tau^{\varepsilon}-t]=p_{\pm}+O(\varepsilon^2),
\end{split}
\end{equation*}
where $p_{\pm}=P(X_{\tau^{\varepsilon}}=x\pm \varepsilon)$ and the last inequality is due to the assumption that $c(\cdot)$ is Lipschitz continuous. For the last term, 
\begin{equation*}
   \begin{split}
\mathbb E_{t,x}[|X_{\tau^{\varepsilon}}-c(\tau^{\varepsilon})|1_{\{\tau^{\varepsilon}=t+\varepsilon\}}]=&\mathbb E_{t,x}[|X_{\tau^{\varepsilon}}-c(t)+c(t)-c(\tau^{\varepsilon})|1_{\{\tau^{\varepsilon}=t+\varepsilon\}}]\\
\le&C\varepsilon P(\tau^{\varepsilon}=t+\varepsilon)=C\varepsilon P(\tilde \tau^{\varepsilon}>t+\varepsilon)=O(\varepsilon^2).
   \end{split} 
\end{equation*}
Combining these estimations, we get that 
\begin{equation}\label{est3}
\mathbb E[k(\tau^{\varepsilon},X_{\tau^{\varepsilon}})]=(p_{-}+p_{+})\varepsilon+O(\varepsilon^2)=\varepsilon+O(\varepsilon^2),
\end{equation}
where we use \eqref{prob_esti} again to get the last equality.
Hence, \eqref{est1},\eqref{est2} and \eqref{est3} yield that 
\begin{equation*}
\lim_{\varepsilon\rightarrow 0}\frac{(\mathbb E_{t,x}[l^c_{\tau^{\varepsilon}}-l^c_t])^2}{\mathbb E_{t,x}[\tau^{\varepsilon}-t]}=\sigma^2(t,x).
\end{equation*}
\section{A Key Lemma}
\begin{lemma}\label{append_lemma}
Let $\varphi$ be the solution of 
$$
(\partial_t+\mathcal L)\varphi+c\varphi=g, \varphi(T,x)=f,
$$
with $c,g$ being bounded continuous functions and $f$ also bounded with bounded derivatives. Then
\begin{enumerate}
    \item There exists a constant $C$ depending on the coefficients such that 
    $$
    -C(\|f^-\|_\infty+T\|g^-\|_\infty)\le \varphi \le C(\|f^+\|_\infty+T\|g^+\|_\infty),
    $$
    where $f^{\pm}$ represents the positive and negative part of $f$ respectively. If $c\le 0$, the constant $C$ can be chosen to be $1$.
    \item Assuming that $c\equiv 0$, we have 
    $$
    \|\partial_x \varphi\|_\infty \le (1+C\sqrt{T})\|\partial_x f\|_\infty+C(\sqrt{T}+T)\|g\|_\infty.
    $$
\end{enumerate}
\end{lemma}
\begin{proof}
From Feymann-Kac representation, it holds that 
$$
\varphi(t,x)=\mathbb E_{t,x}\left[e^{\int_t^T c(u,X_u)du}f(X_T) +\int_t^Te^{\int_t^s c(u,X_u)du}g(s,X_s)ds \right].
$$
Then, one can easily get the first estimation from the assumptions on the coefficients. 

To prove the second estimation, define two processes $\nabla X$ and $N$ as
$$
d\nabla X=\partial_xb(X_t)\nabla X_tdt+\partial_x\sigma(X_t) \nabla X_t dW_t,\nabla X_t=I, 
$$
and 
$$
N_s=\frac{1}{s-t}\int_t^s <\sigma^{-1}(X_u)\nabla X_u ,dW_u>.
$$
One can verify that 
$$
\mathbb E_{t,x}[\sup_{t\le s \le T}|\nabla X_s-I|^2]\le CT\text{, and }\mathbb E_{t,x}[|N_s|^2]\le C(\frac{1}{s-t}+1).
$$
Then, using Bismut-Elworthy-Li formula \cite{elworthy1994formulae}, we get that, for any function $\xi$
$$
\partial_x \mathbb E_{t,x}[\xi(X_s)]=\mathbb E_{t,x}[\partial_x\xi(X_s)\nabla X_s]=\mathbb E_{t,x}[\xi(X_s)N_s].
$$
Thus,
$$
\partial_x \varphi(t,x)=\mathbb E_{t,x}[\partial_x f(X_T)\nabla X_T]+\int_t^T\mathbb E_{t,x}[g(s,X_s)N_s] ds.
$$
Clearly,  we have 
$$
|\mathbb E_{t,x}[\partial_x f(X_T)\nabla X_T]| \le  (C\sqrt{T}+1)\|\partial_x f\|_\infty,
$$
and 
\begin{equation*}
	\begin{split}
		&\int_t^T\mathbb E_{t,x}[g(s,X_s)N_s] ds \le \|g\|_\infty \int_t^T (\mathbb E_{t,x}[|N_s|^2])^{\frac12}ds\\
		 \le& C\|g\|_\infty \int_t^T \frac{1}{(s-t)^\frac12}+1ds=C(\sqrt{T-t}+(T-t))\|g\|_\infty.
	\end{split}
\end{equation*}
This  gives the second estimation. 
\end{proof}

\end{document}